\documentclass[12pt]{amsart}
\usepackage{geometry}                
\usepackage{graphicx, color}
\usepackage{amssymb, amsmath}
\usepackage{epstopdf}
\newtheorem{theorem}{Theorem}[section]
\newtheorem{lemma}[theorem]{Lemma}
\newtheorem{corollary}[theorem]{Corollary}
\newtheorem{proposition}[theorem]{Proposition}

\newtheorem{definition}[theorem]{Definition}

\newcommand{\Z}{\mathbb{Z}}

\newcommand{\R}{\mathbb{R}}

\newcommand{\acts}{\curvearrowright}

\newcommand{\Stab}{\mbox{Stab}}

\title{JSJ decompositions of Quadratic Baumslag-Solitar groups}
\author{Juan Alonso}

\begin{document}
\maketitle

\begin{abstract}{Generalized Baumslag-Solitar groups are defined as  fundamental groups
of graphs of groups with infinite cyclic vertex and edge groups. Forester \cite{forester}
proved that in most cases the defining graphs are cyclic JSJ
decompositions, in the sense of Rips and Sela.  Here we extend
Forester's results to graphs of groups with vertex groups that can be
either infinite cyclic or quadratically hanging surface groups.}
\end{abstract}

\section{Introduction}


To understand a group $G$, it is often useful to decompose it as an amalgamated free product or an HNN extension over a subgroup that belongs to a well-understood class of groups, such as trivial groups, finite groups or cyclic groups. More generally, consider all possible factorizations of $G$ as a graph of groups with edge stabilizers in some single class of groups.

It is often possible to show the existence of a single graph of groups decomposition, from which all of these factorizations can be obtained. This is called a {\em JSJ decomposition} of $G$ (over subgroups in the given class), although the notion is imprecise on how the other factorizations of $G$ are to be obtained from the JSJ decomposition. An example is Grushko's theorem, which gives all the maximal decompositions of a finitely generated group over the class of trivial groups (i.e. the free factorizations).


The letters JSJ stand for Jaco, Shalen and Johannson. Their results in \cite{js} and \cite{jo} can be interpreted as the existence of a JSJ decomposition for $3$-manifold groups over subgroups isomorphic to $\Z\times\Z$. It was these works that motivated the study of JSJ decompositions over non-trivial subgroups (i.e. aside form the Grushko decomposition). Various existence theorems were obtained by Kropholler \cite{k}, Rips and Sela  \cite{sela}, \cite{ripsela}, Bowditch \cite{bow}, Dunwoody and Sageev \cite{dsa}, Fujiwara and Papasoglu \cite{fp}, Dunwoody and Swenson \cite{dsw} and Scott and Swarup \cite{ss}. In  \cite{gl0} and \cite{gl1} Guirardel and Levitt propose a precise definition of JSJ decomposition, which is verified by the graphs of groups constructed in most of the mentioned works. The constructions in \cite{dsw} and \cite{ss} constitute other notions of JSJ decomposition, as does the compatibility JSJ decomposition in \cite{gl2}.


In this paper we will focus on the JSJ decomposition due to Rips and Sela \cite{ripsela}, for finitely presented one-ended groups, with infinite cyclic edge stabilizers (stated here as theorem \ref{rips-sela}).


It is not always clear how to recognize whether a given graph of groups is a JSJ decomposition of its fundamental group. In \cite{forester}, Forester studied the {\em Generalized Baumslag Solitar (GBS) groups}, which are defined by graphs of groups whose vertex and edge stabilizers are infinite cyclic. He showed that the defining graph of a GBS group is a JSJ decomposition, under mild hypotheses. 

In this paper we introduce the {\em Quadratic Baumslag Solitar (QBS) groups}.  They are defined by graphs of groups whose edge groups are infinite cyclic, and whose vertex groups can be either infinite cylic or {\em quadratically hanging} surface groups. (For the meaning of quadratically hanging see definition \ref{QH}). It is clear that the GBS groups are a subclass of the QBS groups. We extend Forester result to the class of QBS groups. I.e. we show that the defining graph of a QBS group is a JSJ decomposition, under some conditions. Specifically, the main theorem in this paper is the following.

\begin{theorem} \label{jsj} Let $\Gamma$ be a QBS graph, $G=\pi_1(\Gamma)$. Suppose that $\Gamma$ is reduced and satisfies the following conditions:
\begin{enumerate}
\item Each edge $e$ of $\Gamma$ has labels $m_e^+,m_e^->1$.
\item Each GBS component $\Gamma_i$ of $\Gamma$ is reduced, and $T_{\Gamma_i}$ is not a point or a line.   
\end{enumerate}
Then $\Gamma$ is a Rips-Sela JSJ decomposition for $G$
\end{theorem}

A {\em QBS graph} is the defining graph of a QBS group. See section \ref{jsjintro:elem} for the definition of a reduced graph of groups. For the edge labels, see definition \ref{labs}. And the GBS components of a QBS graph are defined at the begining of section \ref{s-qbs}.  


The paper is organized as follows: In section \ref{jsjintro} we review the basics of the Rips-Sela JSJ decomposition. There we discuss {\em universality} and {\em unfoldedness}, the two main conditions for a graph of groups (with cyclic edge groups) to be a JSJ decomposition (in the sense of Rips and Sela). In section \ref{gbs-section} we recall the results of Forester about GBS groups, which we will need when dealing with QBS groups. In section \ref{s-unf} we prove a general criterion for unfoldedness, theorem \ref{unf-univ}, in the same fashion of Proposition 2.17 in \cite{forester}. Section \ref{s-univ} is devoted to theorem \ref{univ-new}. This theorem applies to general graphs of groups with cyclic edge stabilizers, and it allows us to show the universality of the whole graph from the universality of certain subgraphs. Finally, in section \ref{s-qbs} we introduce the QBS groups, and prove they are one-ended with the exception of $\Z$ (Proposition \ref{one-end}). We also give the proof of theorem \ref{jsj}. In this proof we combine theorem \ref{univ-new} with the results of Forester (section \ref{gbs-section}) to show universality, and we use theorem \ref{unf-univ} to prove unfoldedness. 

\section{Review of the JSJ decomposition} \label{jsjintro}

\subsection{Graphs of groups and Bass-Serre theory}

Bass-Serre theory is the most fundamental prerequisite for understanding the JSJ decomposition. This is a very quick review, for a comprehensive treatment see \cite{serre}. 

For a {\em graph} we understand a pair of sets $\Gamma = (V,E)$, the vertex and edge set of $\Gamma$ respectively, together with two maps $s,t:E\to V$, which give the {\em source} and {\em target} of an edge. If $e\in E$ is such an edge, the vertices $s(e)$, $t(e)$ will be called the {\em endpoints} of $e$. Thus our graphs have oriented edges, and admit loops (edges $e$ with $s(e)=t(e)$) and multiple edges (different edges having the same endpoints). 

\begin{definition} A {\em graph of groups} consists on the following:
\begin{enumerate}
\item A connected graph $\Gamma$.
\item A group $G_v$ for each vertex $v$ of $\Gamma$.
\item A group $G_e$ for each edge $e$ of $\Gamma$, and two injective homomorphisms \[ \partial^+_e:G_e\rightarrow G_{v^+} \] \[ \partial^-_e:G_e\rightarrow G_{v^-} \] where $v^+$, $v^-$ are the endpoints of $e$.
\end{enumerate}
\end{definition}

This is denoted by $(\Gamma,G,\partial^+,\partial^-)$, or simply by $\Gamma$. Note that even if the endpoints of an edge $e$ agree, i.e. $v^+=v^-$, there are two different maps $\partial_e^+$ and $\partial_e^-$, one for the source and one for the target of $e$.

If $T$ is a spanning tree for $\Gamma$, let $\pi_1(\Gamma,T)$ be defined by the following presentation.

\begin{itemize}
\item Generators: the elements of $G_v$ for the vertices $v\in V(\Gamma)$, and an element $t_e$ for each edge $e\in E(\Gamma)$, $e\notin T$.

\item Relations: the relations in $G_v$ for each vertex $v$, and 
\[ \partial_e^+(g)=\partial_e^-(g) \qquad \mbox{for } e\in T,\, g\in G_e \]
\[ t_e\partial_e^+(g) t_e^{-1} = \partial_e^-(g) \qquad \mbox{for } e\in E(\Gamma),\, e\notin T,\, g\in G_e \] 

\end{itemize}

This group is called the {\em fundamental group} of $\Gamma$. It does not depend on the spanning tree $T$:

\begin{proposition} If $T$, $S$ are two spanning trees for $\Gamma$, then $\pi_1(\Gamma,T)\cong\pi_1(\Gamma,S)$.
\end{proposition}

Thus we often drop $T$ from the notation.
When $G$ is a group and $G\cong \pi_1(\Gamma)$ for a graph of groups $\Gamma$, we say that $\Gamma$ is a {\em splitting} of $G$.
Note that one-edge splittings correspond to decompositions of $G$ as an amalgamated product or an HNN extension.

If $\Gamma$ is a graph of groups and $A\subset\Gamma$ is a connected subgraph, let $\Gamma^* = \Gamma/A$ collapsing $A$ to a vertex $w$. Put $G_w = \pi_1(A)$ and leave the same groups in the non-collapsed vertices and edges. This defines a graph of groups in $\Gamma^*$.

\begin{proposition} If $A\subset\Gamma$ is a connected subgraph and $\Gamma^*=\Gamma/A$ as above, then $\pi_1(\Gamma)\cong\pi_1(\Gamma^*)$.
\end{proposition}

We say that $\Gamma$ is a {\em refinement} of the splitting $\Gamma^*$. Through this proposition, we can see general splittings as iteration of amalgamated products and HNN extensions.

The following is the main result in Bass-Serre theory. It relates the splittings of a group with its actions on simplicial trees, and it is crucial for both the theory of JSJ decompositions and the results on this paper.

\begin{theorem}(Bass-Serre) Given a group $G$, there is a correspondence between the splittings of $G$ as a graph of groups and the actions of $G$ on simplicial trees without edge inversions. If $G\acts X$ is such an action, then the corresponding splitting can be constructed as follows:

- The underlying graph is $\Gamma = X/G$. 

- If $\tilde x \in X$ is a vertex or edge, and $x\in\Gamma$ is its projection, then $G_x$ is conjugate to $\Stab_G(\tilde x)$.\\

\end{theorem}

Under this bijection, the action corresponding to a splitting $\Gamma$ will be denoted by $T_{\Gamma}$, the Bass-Serre tree of $\Gamma$.


Given a simplicial $G$-tree $X$, a subgroup $H\leq G$ {\em acts elliptically} on $X$, or is an {\em elliptic subgroup} with respect to $X$, if there is a point in $X$ that is fixed by every element of $H$. Note that if $X=T_{\Gamma}$ is the Bass-Serre tree of a graph of groups $\Gamma$, then a subgroup of $G=\pi_1(\Gamma)$ acts elliptically on $X$ if and only if it is conjugate into one of the vertex groups of $\Gamma$.


Let $X$, $Y$ be simplicial $G$-trees. A {\em morphism} $f:X\to Y$ is a $G$-equivariant map, which can be made simplicial by subdividing the edges of $X$. The following fact is widely known, and not hard to prove.

\begin{proposition} There is a morphism $X\to Y$ if and only if every elliptic subgroup of $X$ is also elliptic in $Y$.
\end{proposition}

\subsection{Elementary deformations, Foldings} \label{jsjintro:elem}

Here we introduce some important transformations on graphs of groups.

Let $\Gamma$ be a graph of groups. Let $e$ be an edge of $\Gamma$ and $v^+$, $v^-$ its endpoints.
First suppose that $v^+ \neq v^-$ and $\partial^-_e$ is an isomorphism. That is, $G_e = G_{v^-} = C$ and $G_{v^+} = A$ with $C\subset A$.
In this situation, the collapse of the edge $e$ is called an {\em elementary collapse}. Note that $v^+$ and $v^-$ are identified to a single vertex $\bar v$, and $G_{\bar v}= A$ (through the isomorphism $A*_C C\cong A$).

The inverse of an elementary collapse is called an {\em elementary expansion}, and these transformations are the {\em elementary deformations}, which were introduced in \cite{forester2}.

We say that $\Gamma$ is {\em reduced} if it does not admit any elementary collapses.

Again, let $e$ be an edge with different endpoints. This time suppose that $G_e = C \subset C_1 \subset A = G_{v^+}$ and $B=G_{v^-}$. Get $\Gamma_1$ from $\Gamma$ by redefining $G_e = C_1$ and $G_{v^-}=C_1*_C B$. We have $\pi_1(\Gamma)=\pi_1(\Gamma_1)$ by the isomorphism $A*_C B \cong A*_{C_1}(C_1*_C B)$. In this case we say that $\Gamma_1$ is a {\em folding} of $\Gamma$, and that the folding occurs at the vertex $v^+$.

There is another case of folding when $e$ is a loop, that is $v^+ = v^- = v$. Let $G_e=C$, $G_v = A$, and suppose that $\partial^+_e(C)\subset C_1 \subset A$. This time make $\Gamma_1$ with $G_v = A*_C t_e C_1t_e^{-1}$ and $G_e = C_1$. The fundamental group is again preserved, and this transformation is also called {\em folding}. Making some abuse of notation, we say that the folding occurs at $v^+$ in the case just described, and at $v^-$ if we use $\partial^-_e$ instead. 

Looking at the Bass-Serre trees, when there is a folding we have a map $T_{\Gamma} \to T_{\Gamma_1}$, simplicial and equivariant. If $x\in T_{\Gamma}$ is a lift of $v^+$ with stabilizer $gAg^{-1}$, then this map identifies the edges coming from $x$ and projecting to $e$, by the action of $gC_1g^{-1}$. Locally at $x$ it looks like ``folding". In \cite{bestv} this is explained from the viewpoint of graphs of groups.

If $e$ is an edge of $\Gamma$, let $\Gamma_e$ be the graph of groups obtained by collapsing the components of $\Gamma-e$.

\begin{definition} A splitting $\Gamma$ is {\em unfolded} when either:
\begin{enumerate}
\item $\Gamma$ has only one edge, and there is no folding onto it. That is, there is no $\Gamma_0$ such that $\Gamma$ is obtained as a folding of $\Gamma_0$. 
\item $\Gamma$ has several edges, and $\Gamma_e$ is unfolded for all of them.
\end{enumerate}
\end{definition}

\subsection{$\Z$-Splittings, Quadratically hanging subgroups}

A $\Z$-{\em splitting} of the group $G$ is a splitting whose edge groups are infinite cyclic. That is, a graph of groups $\Gamma$, with $\pi_1(\Gamma)\cong G$ and $G_e\cong\Z$ for all edges of $\Gamma$.

\begin{definition} \label{QH} Let $\Gamma$ be a graph of groups. A vertex group $G_v$ is {\em quadratically hanging} (QH) if
\begin{enumerate}
\item $G_v \cong \pi_1(S)$ where $S$ is a 2-orbifold. That is to say, it has one of the following presentations
\[ \langle a_1,\ldots,a_g,b_1,\ldots,b_g,p_1,\ldots,p_m,s_1,\ldots,s_n | s_i^{k_i}=1, \Pi_k p_k \Pi_i s_i \Pi_j [a_j,b_j] =1 \rangle \]
\[ \langle a_1,\ldots,a_g,p_1,\ldots,p_m,s_1,\ldots,s_n | s_i^{k_i}=1, \Pi_k p_k \Pi_i s_i \Pi_j a_j^2 =1 \rangle \]

\item The edges from $v$ are in correspondence with the components of $\partial S$. Moreover, if these edges are $e_1,\ldots,e_m$, then we have $\partial_{e_i}:G_{e_i} \to \langle p_i \rangle$ (where $p_i$ is the boundary loop corresponding to $e_i$). 
\end{enumerate}
\end{definition}

\begin{definition} Let $G$ be a group. Then $P\subset G$ is a {\em QH subgroup} if there is a $\Z$-splitting $\Gamma_P$ of $G$ with $P$ occuring as a QH vertex group.
\end{definition}

Let $\Gamma_1$, $\Gamma_2$ be one-edged $\Z$-splittings of $G$, with edge groups $C_1$, $C_2$ respectively. That is, $G$ is written as an amalgamation or HNN extension over $C_i$. We say that $\Gamma_1$ is {\em elliptic} in $\Gamma_2$ if the subgroup $C_1$ acts elliptically in $T_{\Gamma_2}$, the Bass-Serre tree of $\Gamma_2$. Otherwise, we say that $\Gamma_1$ is {\em hyperbolic} in $\Gamma_2$.

\begin{proposition} (\cite{ripsela}, Theorem 2.1) Let $G$ be freely indecomposable, and $\Gamma_1$, $\Gamma_2$ be one-edged $\Z$-splittings of $G$.  Then $\Gamma_1$ is elliptic in $\Gamma_2$ if and only if $\Gamma_2$ is elliptic in $\Gamma_1$.
\end{proposition}

\subsection{The Rips-Sela JSJ decomposition}

We will now state the fundamental theorem of Rips and Sela, which proves the existence of certain $\Z$-splittings that will be called JSJ decompositions. It applies to {\em one-ended} groups, that are defined as follows.

\begin{definition} A space $X$ is {\em one-ended} if there is a sequence of compact sets $K_n$, such that $X=\cup_n K_n$ and $X-K_n$ is connected for all $n$.
\end{definition}

\begin{definition} A group $G$ is {\em one-ended} if one/all of its Cayley graphs is/are one-ended. Equivalently, for $G$ finitely generated, if it acts freely and cocompactly in a one-ended space.
\end{definition}

Consider a class of groups $\mathcal A$, such as trivial, finite or cyclic groups. We say that a group $G$ {\em splits over $\mathcal A$} if it admits a non trivial graph of groups decomposition with edge groups in $\mathcal A$. For example, $G$ splits over infinite cyclic groups if it admits a non trivial $\Z$-splitting.

According to a theorem of Stallings \cite{st}, a finitely generated infinite group is one-ended if and only if it does not split over finite groups. Thus it makes sense to study the splittings over infinite cyclic groups, $\Z$-splittings, of such a group as a next step. 

\begin{definition} A simple closed curve in a 2-orbifold $S$ is {\em weakly essential} if it is not nullhomotopic, nor boundary parallel, nor the core of a Moebius band embedded in $S$, and does not circle around a branching point.
\end{definition}

\begin{theorem}(Rips-Sela) \label{rips-sela} Let $G$ be a finitely presented one-ended group. Then there is a reduced, unfolded $\Z$-splitting $\Gamma$ of $G$ satisfying the following conditions:
\begin{enumerate}
\item \label{c1}
\begin{itemize}
\item[(a)] A vertex group of $\Gamma$ can either be a QH vertex group, or be elliptic in every $\Z$-splitting of $G$.
\item[(b)] Edge groups are elliptic in every $\Z$-splitting of $G$. 
\item[(c)] Every maximal QH subgroup of $G$ is conjugate to a QH vertex group of $\Gamma$.
\end{itemize}
\item \label{c2} Let $\Gamma_1$ be a one-edged $\Z$-splitting of $G$, with edge group $C$. Suppose that $\Gamma_1$ is hyperbolic in some other one-edged $\Z$-splitting. Then there is a QH vertex group $G_v=\pi_1(S)$ of $\Gamma$, and a weakly essential simple closed curve $\gamma\subset S$ such that $C$ is conjugate to the group generated by $[\gamma]\in G_v\subset G$.

\item \label{c3} If $\Gamma_1$ is a one-edged $\Z$-splitting of $G$ that is elliptic in every other one-edged $\Z$-splitting, then there is a morphism $T_{\Gamma}\to T_{\Gamma_1}$.

\item \label{c4} Let $\Gamma_1$ be any $\Z$-splitting of $G$. Then there is a $\Z$-splitting $\hat\Gamma$, which is a refinement of $\Gamma$ obtained by splitting some QH vertex groups along weakly essential simple closed curves, and a morphism $T_{\hat\Gamma}\to T_{\Gamma_1}$. 
\end{enumerate} 
\end{theorem}

A splitting $\Gamma$ as in the theorem is called a {\em cyclic JSJ decomposition}, or {\em Rips-Sela JSJ decomposition} of $G$. Here we will consider only this version of JSJ decomposition.

Condition \ref{c4} in the theorem is called {\em universality}. It says how every $\Z$-splitting of a group $G$ can be obtained from a JSJ decomposition. Also, it is because of universality that the splitting in the theorem verifies the general definition of a JSJ decomposition (over infinite cylic groups), given by Guirardel and Levitt in \cite{gl0} and \cite{gl1}. Although we will not need that definition here.

There is some redundancy in the conditions for a Rips-Sela JSJ decomposition, as the following proposition shows.

\begin{proposition} \label{universality0} Let $G$ be a one-ended group. Suppose $\Gamma$ is a reduced $\Z$-splitting of $G$ satisfying universality (condition \ref{c4} of \ref{rips-sela}). Then it also satisfies conditions \ref{c1}, \ref{c2} and \ref{c3} of \ref{rips-sela}.
\end{proposition}

{\em Proof:}


For \ref{c1}(a) and \ref{c1}(b), let $\Gamma_1$ be any $\Z$-splitting of $G$. Let $\hat \Gamma$ and $f:T_{\hat\Gamma}\to T_{\Gamma_1}$ be the refinement and the morphism given by universality. If $G_v$ is a vertex group of $\Gamma$ that is not QH, then it is still elliptic in $\hat\Gamma$, and so it is elliptic in $\Gamma_1$. This proves \ref{c1}(a). The edge groups of $\Gamma$ are also elliptic in $\hat\Gamma$, and so they are elliptic in $\Gamma_1$. This gives \ref{c1}(b).

Now we prove condition \ref{c2}. Let $\Gamma_1$ be a one-edged $\Z$-splitting of $G$ that is hyperbolic in some other $\Z$-splitting. Let $\hat \Gamma$ be the refinement of $\Gamma$ given by condition 4, and $f:T_{\hat \Gamma} \to T_{\Gamma_1}$ the corresponding morphism. Take $e$ an edge in $T_{\Gamma_1}$, let $C=\Stab_G(e)$ be its stabilizer subgroup and $K=f^{-1}(e)$ be its pre-image under $f$. There are two kinds of edges in $\hat \Gamma$: those that were already present in $\Gamma$, and those that were obtained by cutting the surfaces of QH vertices along simple closed curves. Since $f(K)=e$, $K$ is not a single point and it meets the interior of an edge $e_1$. Then $\Stab_G(e_1)\subset C$. Moreover, since $C$ is cyclic, the generator of $\Stab_G(e_1)$ is a power of the one of $C$. If $e_1$ was of the first kind, then $C$ would be elliptic in every $\Z$-splitting of $G$, which is a contradiction against our assumption on $\Gamma_1$. Thus $e_1$ is of the second kind, and $K$ does not meet the interior of any edges of the first kind. Let $K^+$ be the union of the edges $e'$ of $T_{\hat \Gamma}$ so that $\Stab_G(e')$ intersects $C$ in a non-trivial subgroup. Then $K^+$ is connected and contains $K$. (If $C=\langle c \rangle$, then $K^+ = \bigcup_{n\geq 1}\mbox{Fix}(c^n)$ which is an increasing union of connected sets). The same reasoning used for $e_1$ shows that $K^+$ does not contain edges of the first kind. (Recall that an element $g$ is elliptic if and only if $g^n$ is elliptic for any $n\neq 0$).

 Now let $v$ be the QH vertex of $\Gamma$ that corresponds to $e_1$. Let $\Gamma_0$ be the splitting of $G_v = \pi_1(S)$ obtained by cutting $S$ along the same simple closed curves as in $\hat \Gamma$. Then there is a copy of $T_{\Gamma_0}$ embedded in $T_{\hat \Gamma}$ that contains $e_1$. Notice that if $g:T_{\hat \Gamma} \to T_{\Gamma}$ is the map that collapses all edges of the second kind, then $g$ collapses $T_{\Gamma_0}$ to a vertex $w$ in the orbit of $v$. So $\Stab_G(T_{\Gamma_0})=\Stab_G(w)$ and it is conjugate to $G_v = \pi_1(S)$. Observe that  $K^+$ must be contained in $T_{\Gamma_0}$, since it can't cross edges of the first kind. In particular, any fixed point of $C$ lies in $T_{\Gamma_0}$, and so it is mapped to $w$ by $g$. Thus $C\subset \Stab_G(w)$ that is conjugate to $ G_v$. And $C=\Stab_G(e_1)$, since a simple closed curve represents a primitive element of $\pi_1(S)$. This proves condition \ref{c2}.

Also note that $K$ does not meet the interior of any other edge of $T_{\Gamma_0}$, since different edges of such tree have different stabilizers. Thus $K$ is contained in $e_1$. Since non-trivial conjugates of $C$ do not stabilize $e_1$, we have that $f(e_1)$ is contained in $e$. Thus $K=e_1$ and the map $f$ results from the collapse of the components of $\hat \Gamma$ minus the interior of the projection of $e_1$. 

Now lets prove condition \ref{c3}. The setup is the same as in the previous case, but this time $\Gamma_1$ is elliptic in every $\Z$-splitting of $G$. If $K$ intersects the interior of an edge of the second kind, we can reason as before and obtain that $C$ is generated by $[\alpha] \in \pi_1(S) = G_v$, a simple closed curve, not in the boundary of $S$. Let $\beta$ be a simple closed curve in $S$ that intersects $\alpha$ non-trivially and minimally. Then the one-edged splitting $\Gamma_2$ of $G$ obtained from $[\beta]$ is hyperbolic in $\Gamma_1$, against our assumption. So $K$ does not intersect any edges of the second kind. This holds for $K=f^{-1}(e)$ where $e$ is any edge in $T_{\Gamma_1}$, so all the edges of the second kind are collapsed to points under $f$. Let $g:T_{\hat \Gamma} \to T_{\Gamma}$ be the map obtained by collapsing the edges of the second kind. Then $f$ factors through $g$, and so we obtain the morphism in condition \ref{c3}.

Finally, for condition 1(c), let $H$ be a QH subgroup of $G$. Let $\Gamma_1$ be a $\Z$-splitting realizing it as a QH vertex. Write $H = \pi_1(S)$ as given by $\Gamma_1$. Again, condition 4 gives a morphism $f:T_{\hat \Gamma} \to T_{\Gamma_1}$ for some refinement $\hat \Gamma$ of $\Gamma$ as before. 

If $c$ is the class of a boundary component of $S$, then $c$ acts elliptically on $T_{\hat \Gamma}$. This is because some power of $c$ fixes an edge $e$ of $T_{\Gamma_1}$ (the incident edge at $v$ corresponding to this boundary curve), and $f^{-1}(e)$ meets the interior of an edge $e_1$ whose stabilizer is also a power of $c$. (By the same argument used to prove condition \ref{c2}).

Consider the action of $H$ on $T_{\hat \Gamma}$ by restriction, and let $\hat T$ be a minimal subtree for this action. Then the boundary classes of $S$ are elliptic in $\hat T$, since they are elliptic in $T_{\hat \Gamma}$.

Consider the decomposition $\Gamma_H$ of $H$ induced by $\hat T$. If $e$ is an edge in $\hat T$, then $\Stab_H(e)\subset \Stab_G(e)$, so the edge groups of $\Gamma_H$ are either trivial or infinite cyclic. Since the boundary classes of $S$ are elliptic in $\hat T$, then $\Gamma_H$ can be extended to $\Gamma_2$, a splitting of $G$ obtained by refining $\Gamma_1$. And since $G$ is one-ended, all edge groups of $\Gamma_2$ are infinite cyclic. Hence all edge groups of $\Gamma_H$ are infinite cyclic.

Using corollary \ref{zvc} (below), $\Gamma_H$ is obtained by splitting $S$ along some disjoint, weakly essential simple closed curves. Now, if $e$ is an edge in $\hat T$, then $\Stab_H(e)$ is generated by a conjugate of one of these curves. So $\Stab_H(e)=\Stab_G(e)$ since the generator of $\Stab_H(e)$ is primitive. And it is also hyperbolic in some $\Z$-splitting of $G$, so $e$ is of the second kind.

We conclude as in the proof of condition \ref{c2}, obtaining that $H$ is conjugate into $G_v$, for $v$ a QH vertex of $\Gamma$. $\Box$

\begin{corollary} \label{universality} Let $G$ be a one-ended group. If $\Gamma$ is a reduced, unfolded $\Z$-splitting of $G$ that verifies universality (condition \ref{c4} from theorem \ref{rips-sela}), then it is a Rips-Sela JSJ decomposition for $G$. 
\end{corollary}

\section{Generalized Baumslag-Solitar groups} \label{gbs-section}

Here we discuss the results in \cite{forester} that are relevant to this paper.


\begin{definition} A {\em Generalized Baumslag-Solitar (GBS) graph} is a graph of groups in which all vertex and edge groups are infinite cyclic.
\end{definition}

Note this is a special case of $\Z$-splitting. A {\em GBS group} is a group obtained as a fundamental group of a GBS graph, and a {\em GBS tree} is the associated Bass-Serre tree. 

\begin{lemma} \label{forester-prep} (Lemma 2.6 in \cite{forester}) Let $\Gamma$ be a GBS graph, $G=\pi_1(\Gamma)$. Assume $G\ncong\Z$, and let $T=T_{\Gamma}$ be the Bass-Serre tree of $\Gamma$. Then:
\begin{enumerate}
\item $G$ is not free.
\item $G$ acts freely on $T\times\R$.
\item $G$ is torsion-free, one-ended and has cohomological dimension 2.
\item $T$ contains an invariant line if and only if $G\cong \Z \rtimes \Z$ (i.e. either $\Z^2$ or the Klein bottle group).
\end{enumerate} 
\end{lemma}

The following is the most general statement about JSJ decompositions of GBS groups.

\begin{theorem} \label{forester-jsj} (Theorem 2.15 in \cite{forester}) Let $\Gamma$ be a GBS graph, $G=\pi_1(\Gamma)$. Suppose $\Gamma$ is reduced, unfolded, and $T_{\Gamma}$ is not a point or a line ($G\ncong \Z,\Z\rtimes\Z$). Then $\Gamma$ is a JSJ decomposition of $G$.
\end{theorem}

In general, it is hard to check wether a splitting is unfolded or not. The following result proves unfoldedness for most GBS graphs.

\begin{proposition} \label{forester-unf} (Proposition 2.17 in \cite{forester}) Let $\Gamma$ be a GBS graph. If every edge group is a proper subgroup of its neighboring vertex groups, then $\Gamma$ is unfolded.
\end{proposition}

The combination the two last statements permits us to recognize most GBS graphs as JSJ decompositions of their fundamental groups.

\section{Criterion for unfoldedness} \label{s-unf}

Here we give a criterion for the unfoldedness of a general $\Z$-splitting. It is a generalization of \ref{forester-unf}, due to Forester, and the proof follows the same lines.

\begin{lemma} \label{unf0} Let $G$ be a freely indecomposable group. Suppose that $\Gamma$ is a $\Z$-splitting of $G$, $e$ is an edge of $T_{\Gamma}$ with endpoints $v_0$, $v_1$ and $H \leq \Stab_G(v_1)$ contains $\Stab_G(e)$ properly. If $\Gamma_1$ is a non trivial unfolding of $\Gamma_e$ at the endpoint $v_0$ of $e$, then $H$ cannot be elliptic in $\Gamma_1$.
\end{lemma}

In the statement of the lemma, we abused notation and still called $e$, $v_0$ and $v_1$ to their respective projections in $\Gamma$ and $\Gamma_e$. Recall that $\Gamma_e$ is the graph obtained from $\Gamma$ by collapsing all edges but the projection of $e$.\\

{\em Proof:} 

Let $X$ be the Bass-Serre tree corresponding to $\Gamma_e$ and $Y$ the one coresponding to $\Gamma_1$. Notice that $X$ can be obtained from $T_{\Gamma}$ by collapsing the components of $T_{\Gamma} - Ge$.

Let $q:T_{\Gamma}\to X$ be the quotient map, and $f:Y \to X$ be the folding map. Let $e'$ be an edge of $Y$, with endpoints $v_0'$, $v_1'$, such that $f(e')=q(e)$ and the fold occurs at $v_0'$.

Let $g$ be the generator of $\Stab_G(e)$ and $g^m$ the one of $\Stab_G(e')$. We know $m\neq 0$ since $G$ is freely indecomposable, and so $|m|>1$ since the fold is non trivial ($\Stab_G(e')$ is strictly contained in $\Stab_G(e)=\Stab_G(q(e))$). We may assume $m>1$, the case for $m<-1$ being analogous.

Define $Y_0$, $g^kY_1$ for $k=0,\ldots,m-1$ to be the components of $Y$ minus the edges $g^ke'$, containing $v_0'$, $g^kv_1'$ respectively. Also let $X_0$, $X_1$ be the components of $X-q(e)$ containing $q(v_0)$, $q(v_1)$,  and $T_0$, $T_1$ the ones of $T_{\Gamma}-e$ containing $v_0$, $v_1$. Observe that $f(Y_0)=X_0$, $f(g^kY_1)=X_1$, $q(T_0)=X_0$ and $q(T_1)=X_1$.

Seeking a proof by contradiction, suppose that $H$ is elliptic in $\Gamma_1$. Thus $H$ fixes a point $x'$ in $Y$. Since $g\in H$, and $g$ fixes no point of $g^kY_1$ for any $k$, we get that $x'$ must belong to $Y_0$. Then $H$ fixes the point $x=f(x')$ in $X_0$, and stabilizes the subtree $q^{-1}(x)$ in $T_0$.

Now, $e$ separates $q^{-1}(x)$ from $v_1$, and $H$ stabilizes both. So $H$ must also stabilize $e$, which is a contradiction, since $H$ contained $\Stab_G(e)$ strictly. $\Box$ \\

We can obtain \ref{forester-unf} from this lemma as follows. \\

{\em Proof of \ref{forester-unf}:}

Suppose $\Gamma$ is a GBS graph in the conditions of \ref{forester-unf}. Notice that if $\Gamma$ is not a single vertex, then $G=\pi_1(\Gamma)\ncong\Z$ and so it is one-ended. If $\Gamma$ is not unfolded, then there is an edge $e$ of $\Gamma$ and a non trivial unfolding $\Gamma_1$ of $\Gamma_e$. In the Bass-Serre tree $T_{\Gamma}$, let $v_0$ be the endpoint of $e$ at which the unfolding occurs, and $v_1$ be the other endpoint. Let $e'$ be the edge of $T_{\Gamma_1}$ with stabilizer contained in $\Stab_G(e)$. Put $H=\Stab_G(v_1)$. Then $\Stab_G(e')\leq\Stab_G(e)\leq H$, where both inclussions are strict (the first one because the unfolding is non trivial, the second one by the hypothesis of \ref{forester-unf}). These three subgroups are infinite cyclic, and $\Stab_G(e')$ is elliptic in $\Gamma_1$, so $H$ must also be elliptic in $\Gamma_1$ (if $g^n$ acts elliptically on a tree, so acts $g$). This contradicts lemma \ref{unf0}. $\Box$  \\ 

The following is the main result of this section. It gives an unfoldedness criterion for universal $\Z$-splittings.

\begin{theorem} \label{unf-univ} Let $G$ be a one-ended group. Suppose that $\Gamma$ is a reduced $\Z$-splitting of $G$ satisfying universality. If every edge group is a proper subgroup of its neighboring vertex groups, then it is unfolded, and is therefore a cyclic JSJ decomposition for $G$.
\end{theorem}

{\em Proof:}

Again, suppose that $\Gamma$ is not unfolded. Let $e$ be an edge of $\Gamma$ and $\Gamma_1$ a non trivial unfolding of $\Gamma_e$. Let $v_0$ and $v_1$ be the endpoints of $e$, when considered in $T_{\Gamma}$, and assume the unfolding occurs at $v_0$.

By the universality of $\Gamma$, it has a refinement $\hat \Gamma$, obtained as in condition \ref{c4} of \ref{rips-sela}, that admits a morphism $T_{\hat\Gamma} \to T_{\Gamma_1}$. Let $w_0$, $w_1$ be the vertices of $e$ as an edge of $T_{\hat \Gamma}$, that correspond to $v_0$, $v_1$ respectively. Put $H=\Stab_G(w_1)$. 

Since $H$ is elliptic in $\hat \Gamma$ and there is a morphism $T_{\hat\Gamma} \to T_{\Gamma_1}$, then $H$ must also be elliptic in $\Gamma_1$.

On the other hand, $H\leq \Stab_G(v_1)$ and it contains $\Stab_G(e)$. If $v_1$ is not a QH vertex, then it doesn't get split in the refinement $\hat\Gamma$. So $H=\Stab_G(v_1)$, which contains $\Stab_G(e)$ strictly by hypothesis. And if $v_1$ is a QH vertex, with $G_{v_1}=\pi_1(S)$, then $H$ is conjugate to $\pi_1(S_0)$ where $S_0$ is a component of $S$ cut by some weakly essential simple closed curves. Thus $H$ is not cyclic, and therefore must contain $\Stab_G(e)$ strictly.

By the lemma, $H$ cannot be elliptic in $\Gamma_1$, which is a contradiction. $\Box$

\section{Adding surface vertices to universal graphs} \label{s-univ}

In this section we deduce the universality of a $\Z$-splitting, given the universality of certain subgraphs of it. We start with some preliminaries.

\begin{definition} \label{labs} Let $\Gamma$ be a $\Z$-splitting of a finitely generated group, and $e$ an edge in $\Gamma$. Let $v^+$ and $v^-$ be the endpoints of $e$, and $a$ be a generator of $G_e$. Define $m_e^+$ as the supremum of the $m$ such that $\partial_e^+(a)=b^m$ for some $b\in G_{v^+}$. Define $m^-_e$ in the same manner. 
\end{definition}

The number $m^+_e$ will be called the {\em label} of $e$ at the endpoint $v^+$. (With some abuse of notation, for when $e$ is a loop, it gets two labels, one for each boundary map). We remark that it is possible to have $m^+_e=+\infty$, although this will not happen in the cases that concern us. If $v^+$ is a QH vertex with $G_{v^+}=\pi_1(S)$, then the element $b$ in the definition is the class of the boundary component of $S$ corresponding to $\partial_e^+$. In particular $m^+_e$ is finite. Also, in the case when $G_{v^+}$ is cyclic, the element $b$ is a generator of $G_{v^+}$ and the label is also finite.  

The following theorem, due to Zieschang, will be crucial in the proof of \ref{univ-new}. The proof is referred, and the corollary results from iterated use of the theorem.

\begin{theorem} \label{zvc0} (Thm 4.12.1 in \cite{zieschang}, pag 140) Let $S$ be a 2-orbifold with boundary components $\gamma_1,\ldots,\gamma_n$. Let $\Delta$ be a one-edged $\Z$-splitting of $\pi_1(S)$ in which $[\gamma_1],\ldots,[\gamma_n]$ are elliptic. Then there is a weakly essential simple closed curve $c$ in $S$ such that $\Delta$ is obtained by cutting $S$ along $c$ (via Van-Kampen's theorem). 
\end{theorem}

\begin{corollary} \label{zvc} Let $S$ be a 2-orbifold with boundary components $\gamma_1,\ldots,\gamma_n$. If $\Delta$ is a general $\Z$-splitting in which $[\gamma_1],\ldots,[\gamma_n]$ are elliptic, then $\Delta$ is obtained by cutting $S$ along $c_1,\ldots,c_m$, disjoint weakly essential simple closed curves.
\end{corollary}

We will also need the following simple lemma about coverings of surfaces and $2$-orbifolds.

\begin{lemma} \label{4cover} Let $S$ be a $2$-orbifold with boundary that is not the disk nor the projective plane. Then there is a $4$-sheeted cover $\hat S$ of $S$, such that every boundary component $\gamma$ of $S$ is covered by two boundary components $\hat \gamma_0$, $\hat \gamma_1$ of $\hat S$, and each one is a double cover of $\gamma$. 
\end{lemma} 

{\em Proof:}

Assume $S$ is an orientable surface, the general case is analogous. Write 
\[ \pi_1(S) =  \langle a_1,\ldots,a_g,b_1,\ldots,b_g,p_1,\ldots,p_m | \Pi_k p_k \Pi_j [a_j,b_j] =1 \rangle \]
Observe that if the genus is positive, then the kernel of the map $\pi_1(S) \to \Z_2$ sending $a_1$ to $1$ and all other generators to $0$ defines a double cover $S_0$ of $S$ where each boundary component of $S$ is covered by two homeomorphic copies of itself.
On the other hand, when $m$ is even, the map $\pi_1(S) \to \Z_2$ sending all $p_i$ to $1$ and $a_j$, $b_j$ to $0$ is a well defined homomorphism, and its kernel gives a double cover $S_1$ of $S$ in which each boundary component of $S$ is covered twice by a single boundary curve of $S_1$. Notice that $S_1$ always has positive genus, by the Euler characteristic computation for a finite cover.

Combining these covers produces the desired $4$-sheeted cover in the cases when $m$ is even, or $m$ is odd but $S$ has positive genus. It remains the case of a sphere with an odd number of punctures. In this case, we have 
\[ \pi_1(S) = \langle p_1,\ldots,p_m | p_1\cdots p_m = 1 \rangle \]
Consider the map  $\pi_1(S) \to \Z_2 \times \Z_2$ that sends $p_1,\ldots,p_{m-2}$ to $(1,0)$, $p_{m-1}$ to $(0,1)$ and $p_m$ to $(1,1)$. This is a well defined homomorphism, and gives the desired covering.

$\Box$

The following result is the main point of this section. Under some conditions, it allows us to recognize the universality of a $\Z$-splitting built from the union of smaller universal graphs and some extra QH vertices.

\begin{theorem} \label{univ-new} Let $\Gamma$ be a $\Z$-splitting of the one-ended group $G$. Let $V=\{v_1,\ldots,v_m \}$ be a subset of the QH vertices of $\Gamma$, and $\Gamma_1,\ldots,\Gamma_k$ be the components of the subgraph spanned by the vertices not in $V$. Put $G_i=\pi_1(\Gamma_i)$. Assume that
\begin{enumerate}
\item \label{con1} If $e$ is an edge with endpoints in $V$, then $m_e^+,m_e^->1$.
\item \label{con2} If the vertices $v_j\in V$ and $w\in\Gamma_i$ are connected by an edge, then $w$ is not a QH vertex of $\Gamma_i$.
\item \label{con3} Each $G_i$ is one-ended, and each $\Gamma_i$ satisfies universality as a $\Z$-splitting of $G_i$.
\end{enumerate}
Then $\Gamma$ satisfies universality.
\end{theorem}

{\em Proof:}

First we observe that if $w$ is a vertex of $\Gamma_i$, then it is QH in $\Gamma$ if and only if it is QH in $\Gamma_i$. If $w$ is QH in $\Gamma_i$, then it has no more incident edges in $\Gamma$ by condition \ref{con2}, and so it is also QH in $\Gamma$. And if $w$ is QH in $\Gamma$, then it cannot be connected by an edge to $v_j\in V$, for that would cause $G_i$ to be freely decomposable.

Let $\Gamma'$ be a $\Z$-splitting of $G$, and $T' = T_{\Gamma'}$ its Bass-Serre tree. Consider the action of $G_i$ on $T'$ by restriction of the action of $G$. Passing to a minimal invariant subtree, $G_i$ acts cocompactly and with cyclic edge stabilizers (since $G_i$ is one-ended). So this action gives rise to a $\Z$-splitting of $G_i$. By universality of $\Gamma_i$, there is a refinement $\hat \Gamma_i$, and a morphism $T_{\hat \Gamma_i} \to T'$, so that $\hat \Gamma_i$ is obtained from $\Gamma_i$ by splitting QH vertex groups along weakly essential simple closed curves. Then all the non-QH vertex groups, and all the edge groups of $\Gamma_i$ are elliptic in $\Gamma'$.

This proves that all the non-QH vertex groups of $\Gamma$ are elliptic in $\Gamma'$, since $V$ consists only of QH vertices. 

It also implies that if an edge $e$ has an endpoint in some $\Gamma_i$, then $G_e$ is elliptic in $\Gamma'$: If $e$ is contained in $\Gamma_i$ we have already shown it. If $e$ has endpoints $v_j\in V$ and $w\in\Gamma_i$, then $w$ is non-QH by condition \ref{con2}, and so $G_w$ is elliptic in $\Gamma'$. Since $G_e\subset G_w$, then $G_e$ must also be elliptic in $\Gamma'$.\\

{\bf Claim:} All edge groups of $\Gamma$ are elliptic in $\Gamma'$.\\

{\em Proof:}

If $e$ has an endpoint in some $\Gamma_i$, we have already proved it.

Now let $e$ be an edge with endpoints $v^{\pm}\in V$ (which can be the same vertex). 

Let $\gamma^{\pm}$ be the boundary components of the orbifolds $S^{\pm}$ corresponding to $G_{v^{\pm}}$, so that $\partial^{\pm}_e:G_e \to \langle [\gamma^{\pm}]\rangle$. Let $H_e = \langle [\gamma^+],[\gamma^-] \rangle \subset G$ be the subgroup generated by the classes of $\gamma^{\pm}$. Note that $H_e$ is a GBS group.

If either $m_e^+>2$ or $m_e^->2$, then the splitting of $H_e$ with edge $e$ satisfies the conditions in Forester's theorem (\ref{forester-jsj}), that are direct consequences of those over $m_e^{\pm}$. So it is a JSJ decomposition of $H_e$, and so $G_e$ is elliptic in $T'$ (as we have done for the $\Gamma_i$). 

If $m_e^+=m_e^-=2$, we proceed by contradiction. Suppose $G_e$ is hyperbolic in $\Gamma'$, and let $c$ be the generator of $G_e$. Take an edge $e'$ of $\Gamma'$ that has a lift to $T'$ lying on the axis of $c$. Then $\Gamma''=\Gamma'_{e'}$ is a one-edged $\Z$-splitting of $G$ in which $c$ is hyperbolic. Let $T''=T_{\Gamma''}$ be its Bass-Serre tree, and let $a$ be the generator of the edge group of $\Gamma''$.

On one hand, we consider the subgroup $H_e$. Note that $H_e=\pi_1(K)$ where $K$ is a Klein bottle. ($K$ is obtained by gluing two M\"obius bands by their boundaries. In this case $\gamma^+$ and $\gamma^-$ are the core circles of the M\"obius bands, and $c$ is their common boundary circle). The action of $H_e$ on $T''$ by restriction gives rise to a $\Z$-splitting of $H_e$ (for $H_e$ is freely indecomposable). Note that $c$ is hyperbolic in it, since it is so in $T''$. So this $\Z$-splitting is non trivial, and we can take $b\in H_e$ a generator of an edge group. Now observe that the edge groups of this decomposition of $H_e$ are all conjugate in $G$ into $\langle a\rangle$. This is so because the only elements that fix an edge of $T''$ are the conjugates of a power of $a$. So we obtain an element $b\in H_e$, $b\neq 1$, which is conjugate to a power of $a$.

On the other hand, we consider the subgroup $M$ constructed as follows. 

\begin{figure}[htbp]
\begin{center}
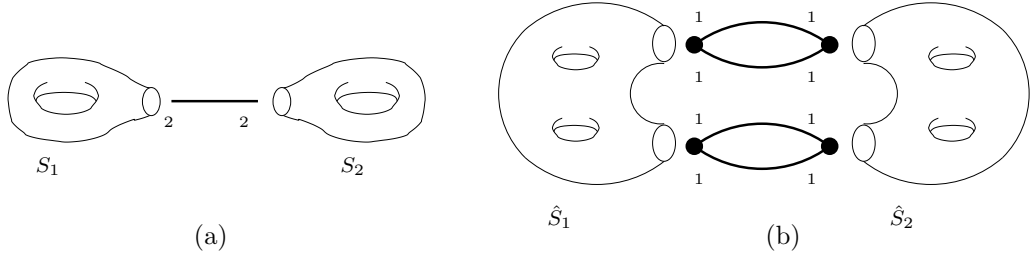
\caption{Example of the cover $\hat \Delta$ {\bf (a)} Original graph $\Delta$, with two QH vertices and an edge with $m_e^{\pm}=2$.  {\bf (b)} Its cover $\hat \Delta$ }
\label{fig0}
\end{center}
\end{figure}

Take the graph formed by the vertices in $V$ and the edges of $\Gamma$ with endpoints in $V$ and both labels equal to $2$. Let $\Delta$ be the component of this graph that contains $e$. For each vertex $v_j\in\Delta$ write $G_{v_j} = \pi_1(S_j)$, where $S_j$ is the orbifold that corresponds to $v_j$ as a QH vertex of $\Gamma$. Let $\hat S_j$ be the $4$-sheeted cover of $S_j$ given by lemma \ref{4cover}. These covers can be extended to a $4$-sheeted cover of the whole graph $\Delta$, that can be constructed as follows. Define the graph $\hat \Delta$ to have the same vertices as $\Delta$, with $\pi_1(\hat S_j) < \pi_1(S_j)$ as vertex group at $v_j$. And for each edge $f$ of $\Delta$, we put in four edges $f_0$, $f_1$, $f_2$ and $f_3$ in $\hat \Delta$, with infinite cyclic edge groups. The boundary maps are described as follows: Suppose $v_j$ is an endpoint of $f$ and $\delta$ is the boundary component of $S_j$ corresponding to $f$. Then let $\delta_0$ and $\delta_1$ be the boundary components of $\hat S_j$ that cover $\delta$ and assign $f_0$, $f_2$ to $\delta_0$ and $f_1$, $f_3$ to $\delta_1$. So the generaror of $G_{f_0}$ maps to $[\delta_0]$ and similarly for the others. This is a $4$-sheeted cover, in the sense that $\pi_1(\hat \Delta) < \pi_1(\Delta)$ with index $4$. (This is best seen by building a presentation 2-complex of $\pi_1(\Delta)$, using $S_j$ for the vertex $v_j$, and tubes for the edges. Then extend the covers $\hat S_j$ of $S_j$ to covers of the tubes.) Note that the labels of the edges of $\hat \Delta$ are all $1$. The local picture at each edge is as in the example on figure \ref{fig0}.

Now let $M$ be the subgroup of $\pi_1(\hat \Delta)$ generated only by the vertex groups and the stable letters of the $f_1$ edges. This is equivalent to saying that $M$ is the fundamental group of the graph resulting from $\hat \Delta$ by erasing all $f_2$ and $f_3$ edges (and keeping the $f_0$ and $f_1$ edges). So $M=\pi_1(S)$, where $S$ is the orbifold that results from gluing the $\hat S_j$ along their boundary curves, so that two boundary curves are identified if they are connected by an edge of $\hat \Delta$. Note that in this subgroup, $c$ is the class of one of the common boundaries of $\hat S^+$ and $\hat S^-$ that corresponds to a lift of the edge $e$. (Say, to $e_0$). Lets call this curve $\beta$, so that $[\beta]=c$.

If $p$ is a boundary curve of $S$, then some power of $[p]$ is in an edge group $G_f$ of $\Gamma$, so that $f$ is not in $\Delta$. (All boundaries corresponding to edges in $\Delta$ were glued). Since $f$ is not in $\Delta$, but connects to a vertex in $\Delta$, we know that $f$ is one of the edges for which we have already proved that $G_f$ is elliptic in $\Gamma'$. Thus, the classes of the boundary curves of $S$ are elliptic in $\Gamma'$. (Ans so in $\Gamma''$).

Again, restrict to $M$ the action on $T''$. This gives a $\Z$-splitting of $M$, in which $c$ is hyperbolic and all the boundary classes of $S$ are elliptic. By corollary \ref{zvc}, this decomposition of $M$ is obtained by cutting $S$ along disjoint, weakly essential simple closed curves. Let $\alpha$ be one of these curves, so that it intersects $\beta$ essentially (i.e. the intersection cannot be removed by homotopy). There must be such $\alpha$, since $c=[\beta]$ is hyperbolic in this decomposition. 

Now, since $[\alpha]$ is a generator of an edge group in the $\Z$-splitting of $M$ induced by $T''$, then $[\alpha]$ must be conjugate in $G$ to a power of $a$. This is by the same argument we used for the element $b$.

Since both $[\alpha]$ and $b$ are conjugate to a power of $a$, then they have the same dynamics in every action of $G$ on a tree. That is to say, in a given $G$-tree, they are either both elliptic or both hyperbolic, depending on the behaviour of $a$.

For the contradiction, consider $\Gamma^*$, the one-edged splitting of $G$ over $[\gamma^-]$. This splitting is obtained from $\Gamma_e$ by folding at $v^+$. In the case of an amalgamation, $\Gamma_e$ corresponds to $A*_{\langle c \rangle}B$ and $\Gamma^*$ to $A*_{\langle[\gamma^-]\rangle}(H_e*_{\langle c\rangle}B)$. The case of an HNN extension is similar.

In both cases $H_e$ is contained in a vertex group, so $b$ must be elliptic in $\Gamma^*$. We will show that $[\alpha]$ is hyperbolic in $\Gamma^*$. This will give the contradiction, thus proving the claim.

Consider the action of $M$ on $T_{\Gamma^*}$ by restriction. It gives a splitting of $M=\pi_1(S)$ in which the boundary classes are elliptic, thus we may use the corollary \ref{zvc} again. This time $c=[\beta]$ stabilizes an edge on $T_{\Gamma^*}$, thus $\beta$ is one of the curves that cut $S$ to form this decomposition. Since $\alpha$ intersects $\beta$ essentially, then $[\alpha]$ must be hyperbolic in this splitting of $M$, and therefore in $\Gamma^*$. $\Diamond$ \\

Thus far we know that all non-QH vertex groups and all edge groups of $\Gamma$ are elliptic in $\Gamma'$. For each QH vertex $v$ of $\Gamma$, write $G_v=\pi_1(S_v)$ where $S_v$ is the corresponding orbifold. Then $G_v$ acts on $T'$ by restriction. Since edge groups of $\Gamma$ are elliptic in $\Gamma'$, it follows that the boundary classes of $S_v$ act elliptically on $T'$. Applying corollary \ref{zvc}, the $\Z$-splitting of $G_v$ induced by its action on $T'$ is obtained by cutting $S_v$ along some disjoint, weakly essential simple closed curves. The vertex groups of this decomposition correspond to the pieces of $S_v$ after the cutting, and are elliptic in $\Gamma'$. Also note that each boundary curve of $S_v$ lies in exactly one of these pieces. So the splitting of $G_v$ is compatible with $\Gamma$, giving rise to a refinement of $\Gamma$. 

Let $\hat \Gamma$ be the refinement of $\Gamma$ that results from splitting all the QH vertex groups $G_v$ as above. Then all vertex and edge groups of $\hat \Gamma$ are elliptic on $\Gamma'$. Equivalently, there is a morphism $T_{\hat \Gamma} \to T'$. Since $\Gamma'$ was an arbitrary $\Z$-splitting of $G$, this concludes the proof. $\Box$

\section{Quadratic Baumslag-Solitar graphs} \label{s-qbs}

Now we consider graphs of groups $\Gamma$ with edge groups infinite cyclic, and vertex groups either QH surface groups or infinite cyclic. We will call these graphs {\em Quadratic Baumslag-Solitar (QBS) graphs}. For simplicity, we restrict the QH vertex groups to be surface groups instead of general 2-orbifold groups. Notice that in a GBS graph all labels are finite, and easily computed from the boundary maps as indicated in the remarks after definition \ref{labs}.

A group $G$ will be called a {\em QBS group} if it can be written as $\pi_1(\Gamma)$, where $\Gamma$ is a QBS graph.

If $\Gamma$ is a QBS graph, let $\Gamma_1,\ldots,\Gamma_k$ be the components of the subgraph spanned by the non-QH vertices. That is, the components that are left after removing all QH vertices and the edges connecting to them. Note that each $\Gamma_i$ is then a GBS graph. The $\Gamma_i$ will be called the {\em GBS components} of $\Gamma$.

\begin{proposition} \label{one-end} Let $\Gamma$ be a reduced QBS graph, and $G=\pi_1(\Gamma)$. Assume $G\ncong\Z$. Then $G$ is one ended.
\end{proposition}

This is a corollary of Theorem 18 in \cite{wilton2}. We also give a proof here.\\

{\em Proof:} Let $X$ be the complex constructed as follows. 

- For each QH vertex $v$ of $\Gamma$, let $X_v=S$ be its corresponding surface, i.e. $G_v=\pi_1(S)$.

- For each cyclic vertex $v$ of $\Gamma$, put in a circle $X_v\cong S^1$.

- For each edge $e$, glue in a cylinder $X_e \cong [-1,1] \times S^1$ along its boundary. The gluing maps are such that they induce $\partial^{\pm}_e$ in the fundamental groups. More explicitly, if $v^{\pm}$ are the endpoints of $e$, then the gluing maps are of the form  $g_e^{\pm}: \{\pm 1 \}\times S^1 \to X_{v^{\pm}}$, so that $(g_e^{\pm})_* = \partial_e^{\pm}$ in the fundamental groups.

By Van-Kampen's theorem, $G=\pi_1(X)$. Let $\tilde X$ be its universal cover. The goal will be to show that $\tilde X$ is one-ended.

Let $\Gamma_1,\ldots,\Gamma_k$ be the GBS components of $\Gamma$. 

Let $X_i \subset X$ be the union of $X_v$, $X_e$ with $v,e\in \Gamma_i$, i.e. the restriction of this complex to the subgraph $\Gamma_i$. 

Notice that the complete lift of $X_i$ consists of disjoint copies of $T_{\Gamma_i}\times \R$ as in lemma \ref{forester-prep}. The complete lift of $X_v\cong S$ for $v$ a QH vertex consists of disjoint copies of $\tilde S$, the universal cover of $S$.

We call the components of the mentioned lifts {\em fundamental pieces}. The fundamental pieces are connected by bands $[-1,1]\times \R$ glued along their boundary lines. In the case of a QH piece $\tilde S$, they are glued to the lifts of the boundaries of $S$. In the case of a GBS piece $\Gamma_i$, they are glued to vertical lines $\{ x\}\times \R$, for $x$ a vertex of $T_{\Gamma_i}$.

Let $\{K_n \}_{n>0}$ be a sequence of compact, simply connected sets that covers $\tilde X$. We need to show that $X-K_n$ is connected for all $n$.
Such sequence exists: a compact set is simply connected if its intersection with each fundamental piece and each band is convex, in the natural geodesic metric of each piece or band.

Let $K\subset \tilde X$ be compact and simply connected. Let $Y$ be a fundamental piece.

If $Y$ is a lift of $X_i$, then $Y\cong T_{\Gamma_i}\times\R$. If $T_{\Gamma_i}$ is not a point, then $Y$ is one-ended, and $Y-K$ has a single component. Thus all points of $Y-K$ are in the same component of $\tilde X - K$. If $T_{\Gamma_i}$ is a point, then there must be a surface (QH) piece connected to $Y$ by a band (if not, $G\cong\Z$). We will deal with this case later on. 

Suppose $Y$ is the universal cover $\tilde S$ of the surface $S=X_v$ for a QH vertex $v$. 

Let $\gamma\subset\partial \tilde S$ be a boundary line. It is connected by a band $B_0$ to another fundamental piece $Y_{\gamma}$. If $Y_{\gamma}$ corresponds to a QH vertex, then $Y_{\gamma} = \tilde S_{\gamma}$, the universal cover of a surface $S_{\gamma}$. Put $B_{\gamma}=B_0$. 

If $Y_{\gamma}$ is a lift of some $X_i$, there must be a QH piece connected to $Y_{\gamma}$ other than $Y$: if not, then $\Gamma_i = \{w \}$ and  $G_e = G_w$ where $e$ is the edge between $v$ and $w$ (so that there are no more lifts of $S=X_v$ connected to $X_i$), and there are no other vertices adjacent to $w$ (no QH pieces projecting to other surfaces). This cannot happen, since $\Gamma$ is reduced. 

Pick one of these QH pieces and call it $\tilde S_{\gamma}$. It is connected by a band $B_1$ to $Y_{\gamma}$. Recall that $Y_{\gamma} \cong T_{\Gamma_i} \times \R$ where $B_0$ is attached to a vertical $\{x_0 \}\times \R$, and so is $B_1$ to $\{x_1 \}\times \R$. Let $\alpha$ be the geodesic of $T_{\Gamma_i}$ between $x_0$ and $x_1$. Then the union of $B_0$, $B_1$ and $\alpha \times \R \subset Y_{\gamma}$ is a band. We call it $B_{\gamma}$.

Thus, for each boundary line $\gamma$ of $Y = \tilde S$, we have another QH piece $\tilde S_{\gamma}$, and a band $B_{\gamma}$ connecting $\tilde S$ to it. Let $D_1$ be the union of all of these pieces and bands. Then $D_1$ is homeomorphic to a disk with some boundary lines (those of the $\tilde S_{\gamma}$ not atteched to $B_{\gamma}$).

We can use the same procedure with the boundary lines of $D_1$ and so on, obtaining $D_k$ for $k\geq 1$. Then $D_k$ is included in the interior of $D_{k+1}$, and each $D_k$ is homeomorphic to a disk with boundary lines.

By compactness, there is some $n$ such that $K \cap D_n$ is in the interior of $D_n$. Since this interior is one-ended, $D_n-K$ has a single component.

So, again, all points of $Y-K$ are in the same component of $\tilde X -K$. The same is true for the points in a band attached to $Y$, thus for the case that was left.  

That covers all the points in $\tilde X -K$, so it has only one component.   $\Box$ \\

This is the main theorem of the paper. It allows us to recognize the defining graph of a QBS group as a Rips-Sela JSJ decomposition, in most cases.

\begin{theorem} (Theorem \ref{jsj}) Let $\Gamma$ be a QBS graph, $G=\pi_1(\Gamma)$. Suppose that $\Gamma$ is reduced and satisfies the following conditions:
\begin{enumerate}
\item \label{cond1} Each edge $e$ of $\Gamma$ has labels $m_e^+,m_e^->1$.
\item \label{cond2} Each GBS component $\Gamma_i$ of $\Gamma$ is reduced, and $T_{\Gamma_i}$ is not a point or a line.   
\end{enumerate}
Then $\Gamma$ is a Rips-Sela JSJ decomposition for $G$
\end{theorem}

{\em Proof:}

Let $V$ be the set of QH vertices of $\Gamma$. The components of $\Gamma$ minus $V$ are the GBS components $\Gamma_i$ of $\Gamma$. By condition \ref{cond1} and \ref{forester-unf}, each $\Gamma_i$ is unfolded. This, together with condition \ref{cond2}, allows us to apply \ref{forester-jsj} (Forester's result). We conclude that each $\Gamma_i$ is a JSJ decomposition of $G_i=\pi_1(\Gamma_i)$. And by \ref{forester-prep}, the $G_i$ are one-ended. By these facts and condition \ref{cond1}, we have verified the hypotheses of theorem \ref{univ-new} for $\Gamma$ and $V$. So $\Gamma$ satisfies universality. Now we can use theorem \ref{unf-univ} to show that $\Gamma$ is unfolded.
Therefore $\Gamma$ is a JSJ decomposition of $G$, by \ref{universality}. $\Box$ \\

\begin{figure}[htbp]
\begin{center}
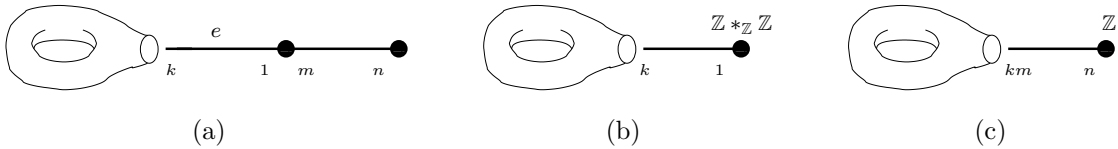
\caption{{\bf (a)} A QBS graph that satisfies universality (for $m,n > 1$), but with $m^-_e = 1$. It admits an unfolding at the surface vertex, as shown in (b) and (c). {\bf (b)} The one-edged splitting corresponding to the edge $e$ of the graph in (a). {\bf (c)} An unfolding of the splitting in (b).}
\label{fig1}
\end{center}
\end{figure}

When some edge label equals $1$, then $\Gamma$ may fail to be a JSJ decomposition. This was already true for GBS graphs. In figure \ref{fig1} there is an example, in which the edge $e$ with a label equal to $1$ is not in a GBS component. However, if in the same figure we change the label $1$ for some $m_e^->1$, and make $k=1$ instead, we do get a JSJ decomposition (by \ref{univ-new} and then \ref{unf-univ}), which is not covered by theorem \ref{jsj}.

\end{document}